\input amstex
 
%
%
%

%
%

%

%

%
\def\zC{{\hbox{$\textstyle C\!\!\!\!I$}}}

%
\def\zN{{\hbox{$\textstyle I\hskip -3pt N$}}}

%

%

%

%

%
%
\def\zR{{\hbox{$\textstyle I\hskip -3pt R$}}}

%

%

%
%
%
\def\zZ{{\hbox{$\textstyle Z\hskip -5pt Z$}}}

%
%
\def\zx{{\raise2pt\hbox{$\chi$}}}
%
%

%
%

%
%
\def\\{$\backslash$}
%
%
\def\subsetnoteq{\subseteq\hskip -14pt\raise-2pt\hbox{
$\scriptscriptstyle\not\phantom{-}$}\hskip 3pt}
%
%
\def\?{\char'76}
%
%

%
%

%
%

%
%
%

%

%
%

%
%

%
%

%
\def\leftitem#1{\item{\hbox to \parindent{\enspace#1\hfill}}}
\def\newline{\relax\ifhmode\null\hfil\break
    \else\nonhmodeerr@\newline\fi}
%
%
\catcode`\@=11
\def\ldisplaylinesno#1{\displ@y\halign{
   \hbox to\displaywidth{$\@lign\hfil\displaystyle##\hfil$}&
     \kern-\displaywidth\rlap{$##$}\kern\displaywidth\crcr
   #1\crcr}}
\catcode`\@=12
%
%
\def\rect{\hskip 0.1truecm\vbox{\hrule\hbox{\vrule\hskip .2truecm
\vbox{\vskip .2truecm}\vrule}
\hrule}\hskip 0.1truecm}
%
%

%
%

%
%
\def\pmb#1{\setbox0=\hbox{#1}%
	   \kern-.025em\copy0\kern-\wd0
	   \kern.05em\copy0\kern-\wd0
	   \kern-.025em\raise.0433em\box 0 }
%
%
%

%
%

\def\hrarr^#1_#2{ \mathrel{ \mathop{\hbox to .5in{\rightarrowfill}} 
 \limits^{\scriptstyle#1}_{\scriptstyle#2}  }}

\def\hlarr^#1_#2{ \mathrel{ \mathop{\hbox to .5in{\leftarrowfill}} 
 \limits^{\scriptstyle#1}_{\scriptstyle#2}  }}

\def\vdarr#1#2{\llap{$\scriptstyle#1$}\left\downarrow\vcenter to .5in{}\right.
\rlap{$\scriptstyle#2$}}



\def\ddrarr#1#2{\llap{$\vcenter {\hbox{$\scriptstyle #1$}}$}
\searrow\displaystyle\rlap{$\vcenter {\hbox{$\scriptstyle #2$}}$}}

%
%
%
%
%
%
\def\bbuildrel#1_#2^#3{\mathrel{\mathop{\kern 0pt#1}\limits_{#2}^{#3}}}
%
%
%

%
%
\def\gordo{\hbox{\bf\kern 0.17em\raise.4ex\hbox{.}\kern 0.17em}}
%
%
\newif\iftitlepage \titlepagetrue
\newtoks\titlepagehead \titlepagehead={\hfil}
\newtoks\titlepagefoot \titlepagefoot={\hfil}

\newtoks\runningauthor \runningauthor={\hfil}
\newtoks\runningtitle \runningtitle={\hfil}

\newtoks\evenpagehead 
\evenpagehead={\hfil\the\runningauthor\hfil}
\newtoks\oddpagehead
\oddpagehead={\hfil\the\runningtitle\hfil}

\newtoks\evenpagefoot \evenpagefoot={\hfil\tenrm\folio\hfil}
\newtoks\oddpagefoot \oddpagefoot={\hfil\tenrm\folio\hfil}

\headline={\iftitlepage\the\titlepagehead
\global\titlepagefalse
\else\ifodd\pageno\the\oddpagehead
\else\the\evenpagehead\fi\fi}

\footline={\iftitlepage\the\titlepagefoot
\global\titlepagefalse
\else\ifodd\pageno\the\oddpagefoot
\else\the\evenpagefoot\fi\fi}

\def\[{\left[}
\def\]{\right]}
\def\({\left(}
\def\){\right)}
\def\l\{{\left\{}
\def\l\|{\left\|}
\def\r\|{\right\|}
\def\r\}{\right\}}
\def\noi{\noindent}
\def\lam{\lambda}
\def\varp{\varphi}
\def\varep{\varepsilon}
\def\gam{\gamma}
\def\Gam{\Gamma}
\def\und{\underbar}
\def\sc{\scriptstyle}
\def\text{\textstyle}
\def\disp{\displaystyle}
\def\widet{\widetilde}
\def\wideh{\widehat}
\def\H{{\cal H}}


\font\tenfrak=eufm10
\font\fivefrak=eufm5
\font\sevenfrak=eufm7
\newfam\frakfam
\textfont\frakfam=\tenfrak
\scriptfont\frakfam=\sevenfrak
\scriptscriptfont\frakfam=\fivefrak
\def\frak{\fam\frakfam\tenfrak}

\documentstyle{amsppt}

\overfullrule=0pt
\normalbaselineskip=15pt
\normalbaselines
\parskip=2pt
\mathsurround=1.8pt

\def\zC{\hbox{\rm C\hskip -5.4pt\vrule
height 8.0pt width 0.4pt depth -0pt\hskip 4.5pt}}
\def\zR{\hskip 2pt\hbox{\rm R\hskip -10.4pt{\rm I}
\hskip 2.5pt}}
\def\rect{\hskip 0.1truecm\vbox{\hrule\hbox{\vrule\hskip .2truecm
\vbox{\vskip .2truecm}\vrule}
\hrule}\hskip 0.1truecm}
\def\noi{\noindent}
\def\UX{\Cal U_{\Cal L_B(X)}}
\def\csta{C$^*$-algebra}
\def\pro{\wp (A,p)}
\def\kp{\Cal K_p (A)}
\def\lp{\Cal L_p (A)}
\def\cu{\Cal U}
\vglue2truecm
\centerline {\bf PROJECTIVE SPACE OF A C$^*$-MODULE  \rm}
\vglue1truecm
\centerline {\smc Esteban Andruchow,  Gustavo Corach and Demetrio Stojanoff
\footnote{Partially supported by UBACYT TW49 and TX92 and ANPCYT 
PICT 97-2259 (Argentina){}}\rm}
\author E. Andruchow G. Corach and D. Stojanoff \endauthor
\title    PROJECTIVE SPACE OF A C$^*$-MODULE  \endtitle
\vskip1truecm

\newdimen
\normalbaselineskip\normalbaselineskip=11pt\normalbaselines
{\narrower \noi \bf Abstract. \sl
Let $X$ be a right Hilbert C$^*$-module over $A$. We study the geometry and the topology of the
projective space $\Cal P(X)$ of $X$, consisting of the orthocomplemented submodules of $X$
which are generated  by a single element. We also study the geometry of the $p$-sphere
$S_p(X)$ and the natural fibration $S_p(X) \to \Cal P(X)$, where $S_p(X)=\{x\in X : <x,x>=p\}$, for $p\in A$ a projection. The projective space and the $p$-sphere are shown to be homogeneous differentiable spaces of the unitary group of the algebra $\Cal L_A(X)$ of adjointable operators of $X$. The homotopy theory of these spaces is examined.
\par}

\newdimen\normalbaselineskip
\normalbaselineskip=15pt
\normalbaselines

\vskip1truecm
\noi \bf Introduction. \rm
Let $A$ be a \csta ,
and $X$ a right Hilbert C$^*$-module over $A$. Denote by $\Cal L_A(X)$ the \csta \ of bounded  adjointable  operators on $X$. In this paper we examine the topology of the set $\Cal P(X)$ of singly generated and
orthocomplemented submodules of $X$, which we call the projective space of $X$. In the classical setting, when $X$ is a finite dimensional vector space over the complex field ($A=\zC$) carrying a positive definite inner product, the topology of the projective space is given by the quotient map
$$
S^1_X \to \Cal P(X),
$$
where $S^1_X$ is the unit sphere of $X$ and $x\in S^1_X$ is mapped to the ray generated by $x$.

To follow this pattern in the C$^*$-module case, one must first recognize which
are the elements of $X$ which generate  orthocomplemented submodules. These turn out to be the $x\in X$ such that there exists $a\in A$ with $xa$ generating the same module as $x$, and $<xa,xa>=p$ a projection in $A$. So a sphere and a map come up, namely the $p$-sphere
$$
S_p(X)=\{ x\in X: <x,x>=p\}
$$
and the map
$$
\rho :S_p(X)\to \Cal P(X) \ , \ \rho (x)=[x]:=\{xa: a\in A\}.
$$
There are many spheres sitting over $\Cal P(X)$.
This fact is to be expected. Suppose that
$[x(t)]$ is a continuous curve in a reasonable topology in  $\Cal P(X)$, with
generators $x(t)$ chosen so that
$<x(t),x(t)>$ are projections,
then the function $t \mapsto <x(t),x(t)>$ should be continuous.
It follows that  elements $[x],[x']$ in $\Cal P(X)$ with non equivalent
projections $<x,x>$ and $<x',x'>$, can not be joined in $\Cal P(X)$
by a continuous  path. On the other hand, one does not need to take into acount all possible
spheres: if $p\sim q$ (Murray-von Neumann equivalence) then $S_p(X)$ and $S_q(X)$ are
isometrically isomorphic and they are mapped onto the same part of $\Cal P(X)$. Indeed, if
$v\in A$ satisfies $v^*v=p$ and $vv^*=q$, then $S_p(X)v^*=S_q(X)$, $S_p(X)=S_q(X)v$,
and thus elements in both spheres generate the same submodules. In other words, any reasonable
topology on $\Cal P(X)$ would define in it at least as many connected components as there are
classes of projections in $A$.

Another point of view would be to consider, instead of orthocomplemented singly generated submodules, the "rank one" projections in $\Cal L_A(X)$, i.e.
projections of the form $\theta_{x,x}$ (where, as is usual notation,
$\theta_{x,y}(z)=<z,x>y$ for $x,y,z \in X$), with the norm topology. It turns out that both topologies, the one induced by the $p$-spheres and this latter
one, coincide. Therefore, in addition to the map $\rho$ (which will be shown to be a fibre bundle) one has the action of the unitary group of $\Cal L_A(X)$: if $\theta_{x, x}$ is a projection and $U$ is a unitary,
$$
U.\theta_{x, x}=U\theta_{x, x}U^*=\theta_{Ux, Ux},
$$
or equivalently, if $[x] \in \Cal P(X)$,
$$
U.[x]=[Ux].
$$

Once these basic facts are established, we proceed to describe the differentiable 
structure of the spheres and the projective space. This is done in section 2. 
In the next sections we consider certain particular examples, and focus on the 
homotopy theory of these spaces.

In section 3 we examine the case $X=H_A$. Using Mingo's theorem [16] one proves that the
spheres $S_p(H_A)$ are contractible. This fact enables one to relate the homotopy groups 
of the connected components of $\Cal P(H_A)$ with those of the unitary groups 
$\Cal U_{pAp}$ of the algebras $pAp$, for $p$ projections in $A$.

In section 4 we study the partial isometries of $\Cal L_A(H_A)=M(A\otimes \Cal K)$. 
Here a clear distinction is drawn between the cases where the initial projection 
$P$ is "compact" (i.e lies in $A\otimes \Cal K$) or not. In the first case, the 
spheres $S_P(M(A\otimes \Cal K))=S_P(A\otimes \Cal K)$ are contractible. 
In the latter case this is far to be the case. For example, as an easy 
consequence of Kasparov's stabilization theorem , $\pi_0$ of the unit sphere 
$S_1(M(A\otimes  \Cal K))$ is a semigroup which naturally parametrizes the 
equivalence classes of countably generated $A$-modules.

In section 5 we examine the von Neumann algebra case, i.e. $A$ is von Neumann and $X$ 
is selfdual. We prove first that the unitary orbit of a projection in a von 
Neumann algebra is simply connected. As a consequence, the $\pi_1$-groups of 
the connected components of $\Cal P(X)$ are trivial. Also several results 
concerning the homotopy groups of the $p$-spheres and the projective spaces are given. 
For example, the Jones index of an inclusion of II$_1$ factors appears as a 
homotopic invariant of the projective space of the $C^*$-module induced by 
the inclusion.

For basic results and notations on $C^*$-modules, we shall follow E. C. Lance's book 
[14] and the papers  [18] and [19] by W. L. Paschke.

\bigskip
\noi \bf 2. Elementary properties of the $p$-spheres  \rm
\medskip
Let $A$ be a \csta \ and $X$ a right C$^*$-module over $A$, which will be supposed to be
full, i.e. $<X,X>$ dense in $A$. As in the introduction, if $p\in A$ is a projection,
$S_p(X)=\{x\in X: <x,x>=p\}$, and the projective space $\Cal P(X)$ is the set
of submodules of $X$ which are generated by a single element, and are also
orthocomplemented, i.e. they are complemented by their orthogonal modules (with
respect to the inner product in $X$).

If $A$ is unital,   $G=G_{A}$ will denote the group of invertibles of $A$ and
$\Cal U=\Cal U_A$
the unitary group of $A$. Let $p=p^2=p^* \in A$ be a projection. Put
$$
\Cal E_p=\{q \in A : q^2=q^*=q, \hbox{ such that there exists } v\in A 
\hbox{ with } v^*v=p , vv^*=q\}.
$$
That is, $\Cal E_p$ is the set of projections of $A$ which are equivalent to $p$. 
There are many papers considering the geometric structure of the space of 
projections of a \csta  \ (see [21], [8], [24]). In [1] it was shown that  
$\Cal E_p$ consists of a union of connected components of the set of all 
projections of $A$.
\medskip
Let us state some general properties of the spheres $S_p(X)$.
\medskip
\noi \smc \bf Lemma 2.1 \rm \sl Let $x\in X$. Then $\theta_{x , x}$ is a 
projection in $\Cal L_A(X)$ if and only if $<x,x>$ is a projection in $A$.

\smallskip \noi Proof. \rm Suppose that $<x,x>=p$ is a projection. Then
$x=xp$, since $<x(1-p), x(1-p)>=(1-p)<x,x>(1-p)=0$ and therefore
$x(1-p)=0$ (this computation can be done in $A_1$ if $A$ has no unit). Then
$$
\theta_{x , x}^2(y)=x<x,x><x,y>=x<x,y>=\theta_{x , x}(y).$$
On the other hand, suppose that $\theta_{x , x}$ is a projection and put 
$<x,x>=r\ge 0$ in $A$. Then $\theta_{x , x}=\theta_{x , x}^2=\theta_{xr , x}$, 
and therefore
$$
r^2= <\theta_{x , x}(x),x>= <\theta_{xr , x}(x),x>=r^3.$$
Since $r\ge 0$ it must be a projection. \rect
\medskip
Note that, therefore, a singly generated module $\{xa: a\in A\}\subset X$
is orthocomplemented if and only if it has a generator $xa_0$ such that 
$<xa_0,xa_0>$ is a projection in $A$. The result above also shows that 
$S_p(X) \subset Xp \subset X$. Therefore it can be regarded as the unit 
sphere of the module $Xp$ over the algebra $pAp$. Unit spheres were studied 
in [1], where it was shown that they enjoy remarkable geometric properties.
\medskip
\noi \smc \bf Remarks 2.2 \rm \sl \rm \rm Let $p\in A$ be a projection 
different from zero.

\item{ 1)} The sphere $S_p(X)$ is a C$^\infty$ complemented submanifold of
$Xp$, and therefore also of $X$ (see [1]).
\item{ 2)} The unitary group $\Cal U_{\Cal L_{pAp}(Xp)}$ of $\Cal L_{pAp}(Xp)$ 
acts on $S_p(X)$ by means of $$U\bullet x=U(x).$$
The action is smooth and locally transitive. More precisely, two elements of 
$S_p(X)$ lying at distance less than $1/2$ are conjugate by this action.
\item {3)} For any  fixed $x_0 \in S_p(X)$, the map
$$
\pi_{x_0} :\Cal U_{\Cal L_{pAp}(Xp)}\to S_p(X) \quad \pi_{x_0} (U)=U(x_0)
$$
is a homogeneous space, with isotropy group $I_{x_0}=\{V \in 
\Cal U_{\Cal L_{pAp}(Xp)}: V(x_0)=x_0 \}$ . In particular, 
it is a principal bundle.
\medskip
\noi \smc \bf Examples 2.3 \rm \sl \rm \rm

\item{1)} If  $X=A$ is a \csta \  one obtains that the right ideal generated
by $b$,
$bA=\{ba: a\in A\}$,
is a complemented submodule of $A$ if and only if it contains a generator
which is a partial isometry.
In other words, $bA$ is orthocomplemented
in $A$ if and only if there exists $bA=pA$ for some projection $p\in A$.
That is, $\Cal P(A)$ identifies with the space of projections of $A$. In this case the $p$ sphere
$S_p(A)$ consists of the partial isometries of $A$ with initial space $p$.
The elements of closed range are related to the differential geometry of the action of the
group of $G_A$ of invertible elements on $A$ by right multiplication. Namely, for $b\in A$ one has the following equivalent
conditions:
\itemitem{1)} The orbit $bG_A=\{bg : g \in G_A\}$
is a  submanifold of $A$
\itemitem{2)} The mapping $\pi_b : G_A \to bG_A$, given by $\pi_b(g)=bg$ is
a C$^\infty$ submersion 
\itemitem{3)} $bA$ is closed in $A$.
\itemitem{4)} $bA=pA$ for some projection $p\in A$
\itemitem{5)} $b$ is relatively regular in $A$.

\item{2)} Suppose that $X=B$ is a unital \csta \ such that $B$ contains $A$ 
and there exists a finite index conditional expectation $E:B\to A$ defining 
the inner product in the usual way ([5], [14], [10]): $<b,b'>=E(b^*b')$. 
In this case 2.1 reads: the submodule 
$\{ba: a\in A\}$ is complemented in $B$ if and only if it contains a generator 
$c$ such that $E(c^*c)$ is a projection in $A$.
The unit sphere $S_1^E(B)=\{x\in B: E(x^*x)=1\}$
(the superscript $E$ stands to distinguish this sphere from the space of
isometries $S_1(B)$, which in fact lies inside $S_1^E(B)$). The final projection $\theta_{1, 1}$
corresponding to the element $1\in S_1^E(B)$ is the Jones projection $e\in
B_1$ of the conditional expectation $E$, where $B_1$ denotes the basic extension
of $A\subset B$ using $E$. 
\medskip

The spheres $S_p(X)$ enable one to endow $\Cal P(X)$ with a natural topology. Namely,
the quotient topology given by the maps
$$
\rho :S_p(X) \to \Cal P(X), \quad \rho (x)=[x]=\{xa: a\in A \}.
$$
Note that if $x,y \in S_p(X)$ are such that $\rho (x)=\rho (y)$, then
$<x,y>$ is a unitary element of $pAp$ verifying $x<x,y>=y$. Indeed,
$y=xa$ for some a in $pAp$. Then $p=<y,y>=a^*<x,y>=<y,x>a$.
Therefore the connected components of $\Cal P (X)$  are identified with the quotient spaces
$$
S_p(X) / \Cal U_{pAp}$$
where $\Cal U_{pAp}$ acts on $S_p(X)$ by the original right $A$-module action.

On the other hand, there is another natural mapping to consider,
$$
\rho': S_p(X) \to \Cal E = \{ \hbox{projections of } \Cal L_A(X)\} , 
\quad \rho'(x)=\theta_{x , x}.
$$
This map is C$^\infty$ and its range is the set of projections which are 
(Murray-von Neumann) equivalent to $\theta_{x_0 ,x_0}$ for any
fixed $x_0 \in S_p(X)$. Denote this set of projections by $\Cal E _{x_0}$. 

\medskip
\noi \smc \bf Proposition 2.4 \rm \sl
If $x_0 \in S_p(X)$, the mapping $\rho' :S_p(X)\to \Cal E_{x_0}$
is a principal bundle with structure group equal to the unitary group of $pAp$.
\smallskip \noi Proof. \rm
It is clear that the range of $\rho'$ is $\Cal E_{x_0}$.
In order to prove that it is a fibre bundle, let us show that $\rho'$ has local 
cross sections around any  $\theta_{x , x}$.
If $P$ is a projection of $\Cal L_A(X)$ with $\|P-\theta_{x , x}\|<1$, then
there exists a unitary operator $U$ (which is an explicit C$^\infty$ formula
of $P$) such that $P=U\theta_{x , x}U^*$ (see [8]). In other words, 
$P=\theta_{U(x) , U(x)}$. Therefore, if $\|\rho'(y)-\rho'(x)\|<1$ , then
$\rho'(y) \mapsto U(x)$ is a local cross section for $\rho'$.
The fibre over $\theta_{x, x}$ equals the set $\{ y \in S_p(X): \theta_{y , y}=
\theta_{x , x}\}$.
Observe that $\theta_{y , y}=\theta_{x , x}$ if and only if 
$y=xu$ for some unitary element $u$ in $pAp$. Indeed, if $y=xu$ it is straightforward to see that $\theta_{y , y}=\theta_{x , x}$. On the other hand, if $\theta_{y , y}=\theta_{x, x}$,
then it is clear that $u=<x,y>$ satifies the required condition. \rect
\medskip
We may state the relation between $\rho$ and $\rho'$ in the following:
\medskip
\noi \smc \bf Theorem 2.5 \rm \sl
Let $X$ be a $C^*$-module over $A$ and $x_0 \in S_p(X)$.  Denote by 
$[x_0]=\{x_0a: a\in A\} \in \Cal P(X)$. 
Then the connected component $\Cal P(X)_{[x_0]}$ of $[x_0]$
in $\Cal P(X)$ is homeomorphic to the connected component of
$\theta_{x_0 , x_0}$ in the space of projections
$\Cal E_{ x_0}$ of $\Cal L_B(X)$.

\smallskip \noi Proof. \rm The proof follows by observing that
both spaces are homeomorphic to the connected component of the
class of $x_0$ in the quotient $S_p(X) /\Cal U_{pAp}$. \rect
\medskip
One may summarize the situation in the following commutative diagram:

\noi
Fix $x_0 \in S_p(X)$, and $[x_0]$ the submodule generated by $x_0$,
$${\normallineskip=6pt\normalbaselineskip=0pt
\matrix
S_p(X)_{x_0}   & \hrarr^{\rho'}_{} & \Cal E_{ x_0}   \\
           & \ddrarr {\rho}{} &  \vdarr {}{r} \\
           &                   &  \Cal P(X)_{[x_0]}
\endmatrix }
$$
where $\Cal P(X)_{[x_0]}$ denotes the connected component of $[x_0]$ in
$\Cal P(X)$, $S_p(X)_{x_0}$ the connected component of $x_0$ in $S_P(X)$,
 and $r$ is a homeomorphism, given by $r(\theta_{x , x})=$range of $\theta_{x , x}$.
for $u \in \Cal U (M)$.
\medskip
One does not need to take into consideration all the spheres $S_p(X)$ for all
possible projections $p \in A$. If $p$ is equivalent to $q$ (for the Murray-von Neumann equivalence), then the spheres $S_p(X)$ and $S_q(X)$ are
diffeomorphic and are carried by $\rho$ onto the same part of $\Cal P(X)$.
Indeed, let $v\in A$ such that $v^*v=q$ and $vv^*=p$. Then, the C$^\infty$ map
$$
S_p(X) \ni x \mapsto xv \in S_q(X)
$$
has inverse
$$
S_q(X) \ni x \mapsto xv^* \in S_p(X),
$$
and clearly, $\rho (x)=\rho(xv)$.
\medskip
Therefore, in order to cover $\Cal P(X)$ one only needs to consider spheres
$S_p(X)$ for projections $p$ chosen one from each equivalence class. The next 
result states that such spheres and their corresponding parts of $\Cal P(X)$ 
lie separated (considering in $\Cal P(X)$ the norm metric of $\Cal L_A(X)$).
\medskip
\noi \smc \bf Proposition 2.6 \rm \sl
Suppose that $p$ and $q$ are projections of $A$ which are not equivalent. Then
\item{1)} $d(S_p(X),S_q(X))\ge \sqrt2 -1$.
\item{2)} $d(\rho'(S_p(X)), \rho'(S_q(X)))\ge 1$.
\smallskip \noi Proof. \rm
a): Pick $x\in S_p(X)$ and $y\in S_q(X)$, and suppose that $\|x-y\|<\sqrt2 -1$. Then
$$
p=<x,x>=<(x-y)+y,(x-y)+y>=\|x-y\|^2+q+<x-y,y>+<y,x-y>.
$$
Therefore
$$
\|p-q\|\le \|x-y\|^2 +2\|x-y\|\|y\|.
$$
Using that $\|y\|=1$ and $\|x-y\|<\sqrt2 -1$, one gets that $\|p-q\|<1$.
This would imply that $p$ and $q$ are unitarily equivalent, a contradiction.

\noi
b): We claim that $\rho'(x)=\theta_{x , x}$ and $\rho'(y)=\theta_{y , y}$
are not unitarily equivalent if $x\in S_p(X)$ and $y\in S_q(X)$ with $p$ and
$q$ non equivalent. Suppose that there exists a unitary $U$ in $\Cal L_A(X)$ such that $U\theta_{x , x}U^*=\theta_{U(x) , U(x)}=\theta_{y , y}$. Then $v=<U(x),y>$ is a partial isometry in $A$ such that $vv^*=<x,x>=p$ and $v^*v=<y,y>=q$.
Indeed,
$$
vv^*=<U(x),y><y,U(x)>=<U(x),\theta_{y , y}(U(x))>=<U(x),\theta_{U(x) , U(x)}(U(x))>=p,$$
and analogously for $v^*v=q$. Note that our claim implies that
$\|\theta_{x, x }-\theta_{y , y } \|\ge 1$. \rect
\medskip
If $p \in A$ is a projection, denote by $\Delta_p(X)$ the set
$$
\Delta_p (X)=\{x\in Xp : <x,x> \hbox{ is invertible in } pAp \}.
$$
\medskip
Clearly  $S_p(X) \subset \Delta_p(X)$. Also it is clear that $\Delta_p(X)$
is an open subset of $X.p$, and therefore a complemented (analytic) submanifold
of $X$.
\medskip
\noi \smc \bf Proposition 2.7 \rm \sl
$S_p(X)$ is a strong deformation retract of $\Delta_p(X)$.
\smallskip \noi Proof. \rm
Put $F_t(x)=x<x,x>^{-t/2}$, for $x \in \Delta_p(X)$ (the inverse
of $<x,x>^{t/2}$ taken in $pAp$). Clearly $F_0=Id_{\Delta_p(X)}$, $F_t|_{S_p(X)}=Id_{S_p(X)}$ and $F_1:\Delta_p(X) \to S_p(X)$ is a retraction. \rect
\medskip
Let $p\ge q$ be projections in $A$. Consider the map
$$
m_{p,q} :S_p(X) \to S_q (X), \quad m_{p,q}(x)=xq.
$$
Since $S_p(X)$ equals the unit sphere $S_1$ of the $pAp$-module $Xp$, we may
restrict our atention to the unital case and the map
$$
m_q=m_{1,q}:S_1(X) \to S_q(X).$$

Note also that the maps $m_{p,q}$ can also be defined on $\Delta_p(X)$ and
$\Delta_q(X)$. Indeed, if $x=xp \in \Delta_p(X)$, then there exists a positive 
number $d$ such that $<x,x>\ge dp$. Therefore $<xq,xq>=q<x,x>q \ge dqpq =dq$, 
i.e. $xq \in \Delta_q(X)$.

The map $m_p$ may not be surjective. Consider for instance $X=B(H)$ for $H$ a 
separable Hilbert space. If $s$ is the unilateral shift, put $p=ss^*$. Then
$s^* \in S_p (B(H))$, but there exists no isometry $v\in S_1(B(H))$ such that
$vp=s^*$.

\noi \smc \bf Remark 2.8 \rm \sl \rm
We have considered the action of the unitary  group of $\Cal L_{pAp}(Xp)$, in order to use the results of [1] which were stated for unit spheres ($p=1$) of Hilbert modules. Nevertheless, it is clear that also the unitary group of $\Cal L_A(X)$ acts on $S_p(X)$. Moreover, the action admits C$^\infty$ local cross sections. Namely, if $x,y \in S_p(X)$ lie at distance less than $1/2$, then
$$
\|\theta_{x , x}-\theta_{y , y}\|\le \|\theta_{x , x}-\theta_{y , x}\|+\|
\theta_{x , y}-\theta_{y , y}\|\le \|x-y\|\|x\| + \|x-y\|\|y\|<1.
$$
This implies that $\theta_{x , x}$ and $\theta_{y , y}$ are unitarily equivalent, and this equivalence can be implemented via a unitary $U$ which is an explicit C$^\infty$ formula of $\theta_{x , x}$ and $\theta_{y , y}$. In other words $\theta_{U(x) , U(x)}=U\theta_{x , x}U^*=\theta_{y , y}$. Then $W=\theta_{U(x) , y}+ 1-\theta_{y , y}$ is a unitary operator of $\Cal L_A(X)$ which verifies $W(y)=U(x)$. Therefore the unitary operator $W^*U$ carries $x$ to $y$, and it is an explicit C$^\infty$ formula of the parameters $x$ and $y$.

In particular, this implies that if $x_0$ is a fixed element of $S_p(X)$, then also the map
$$
\pi_{x_0} : \Cal U_{\Cal L_A(X)} \to S_p(X)_{x_0} \quad \pi_{x_0}(U)=U(x_0)
$$
is a principal bundle, with structure group $\{V\in \Cal U_{\Cal L_A(X)}: V(x_0)=x_0\}$.
\medskip
We may state now the main properties of $m_p$ in the general case.

\medskip
\noi \smc \bf Proposition 2.9 \rm \sl The map $m_p$ fills
connected components, i.e. if $x$ belongs to the image of $m_p$, 
then every other element in the connected component
of $x$ in $S_p(X)$ also belongs to the image of $m_p$.
\smallskip \noi Proof. \rm
Suppose that $y\in S_p(X)$ lies in the same component of $x$.
Put $x_1 \in S_1(X)$ with $m_p(x_1)=x_1p=x$. Since by 2.8 the
action of the unitary group of $\Cal L_A(X)$ on $S_p(X)$ is
locally transitive, there exists a unitary operator $U$ such
that $U(x)=y$. Put $y_1=U(x_1)$. Then clearly $y_1 \in S_1(X)$
and $m_p(y_1)=U(x_1)p=U(x_1p)=y$. \rect
\medskip
\noi \smc \bf Proposition 2.10 \rm \sl
Suppose that $A$ is unital and $m_p(x_1)=x$ for $x\in S_p(X)$
and $x_1\in S_1(X)$. Denote by $S_p(X)_x$ the connected component of $x$
in $S_p(X)$, and analogously for $S_1(X)_{x_1}$. Then the C$^\infty$ map
$$
m_p :S_1(X)_{x_1} \to S_p(X)_x
$$
is a fibre bundle.
\smallskip \noi Proof. \rm
Since the connected component of the identity in the unitary group of $\Cal L_A(X)$ 
acts transitively on both $S_1(X)_{x_1}$ and $S_p(X)_x$, it suffices to show that 
$m_p$ has local cross sections around $x$. If $x' \in S_p(X)$ satisfies 
$\|x-x'\|<1/2$, then by 2.8 there exists a unitary operator $U$ depending 
continuously (in fact, smoothly) on $x'$ such that $U(x)=x'$. Put $s(x')=U(x_1)$. 
Then it is apparent that $s$ is a C$^\infty$ local cross section for $m_p$.\rect
\bigskip

\noi \bf 3. The case $X=H_A$.  \rm
\medskip
Let $H_A$ be the module introduced by Kasparov
$$
H_A=\{ (a_n) : a_n \in A \hbox{ and } \sum_{n \in \zN} a_n^*a_n
\hbox{ converges in norm }\}.
$$
$H_A$ becomes a right $A$-module with the action $(a_n)b=(a_nb)$,
and a Hilbert C$^*$-module with the $A$-valued inner product
$<(a_n),(b_n)>=\sum_{n\in \zN} a_n^*b_n$ (See [13]). The C$^*$-algebra
$\Cal L_A(H_A)$ is isomorphic to $M(A\otimes \Cal K)$, where $\Cal K$
denotes the ideal of compact operators of an infinite dimensional
separable Hilbert space ([13]). We will show that the
$p$-spheres of $H_A$ are contractible, using the Mingo-Cuntz-Higson theorem([16],[9]),
which states that
if $A$ has a countable aproximate unit, then both the invertible and
unitary groups of $\Cal L_A(H_A)$ are contractible. Our result will be valid
though for arbitrary \csta s, and will use only the unital version [16] of the
theorem. We give the statement without proof, since it will be a consequence
of a result in section 4.

\medskip
\noi \smc \bf Theorem 3.1 \rm \sl
For any  \csta $A$ and any projection $p\in A$, $S_p(H_A)$ is contractible. \rect

\medskip
\noi \smc \bf Corollary 3.2 \rm \sl For any projection $p \in A$, the mapping
$m_p:S_1(H_A)\to S_p(H_A)$ is a C$^\infty$ fibre bundle.
\smallskip \noi Proof. \rm It suffices to combine  2.9 and 2.10. \rect
\medskip
\noi \smc \bf Corollary 3.3 \rm \sl For any projection $p\in A$,
the space $\Delta_p(H_A)$ is contractible. \rect \rm
\medskip
As noted before, $S_p(H_A)$ can be regarded as the unit sphere
of the module $H_Ap$. But since clearly $H_Ap$ is different from $H_{pAp}$, the
analysis does not reduce directly to the unital case.

Let us now consider the map $m_p$ on the spaces $\Delta$.
\medskip

\noi \smc \bf Proposition 3.4 \rm \sl
Suppose that $A$ is unital. Then the C$^\infty$ map
$$
m_p :\Delta_1(H_A) \to \Delta_p(H_A)\ , \quad m_p(x)=xp
$$
is a fibre bundle.
\smallskip \noi Proof. \rm
First let us show that $m_p$ is surjective. Let $x=xp \in \Delta_p(H_A)$. Note that if $x=(a_n)$ then $a_np=a_n$ for all $n$. Given $\epsilon >0$, there exists $N$ such that $\|a_N\|<\epsilon$. Pick $x_0=x+e_N(1-p)$ ($e_N$ denotes
the sequence with $1$ in the $N$-th entry and zero elsewhere). Then
clearly $x_0p=x$ and
$$
<x_0,x_0>=1-p+ <x,x>+(1-p)x_N+x_N^*(1-p)
$$
is invertible in $A$, if we choose $\epsilon $ small such that $<x,x>\ge 2\epsilon p$.

\noi
In order to show that $m_p$ is a fibre bundle, it will suffice to show that it has C$^\infty$ local cross sections around any point in $\Delta_p(H_A)$. Note that the connected (in fact, contractible [16]) group
$G_{\Cal L_A(H_A)}$ acts transitively on $\Delta_q$ for any projection $q$.
Indeed, first note that $G_{\Cal L_A(H_A)}$ acts on $\Delta_q$:
if $x\in \Delta_q$ and $G\in G_{\Cal L_A(H_A)}$,
$$
<G(x),G(x)>= <G^*G(x),x> \ge r<x,x>,$$
for some $r>0$, because $G^*G$ is positive and invertible (see [14]).
Next we show that $G_{\Cal L_A(H_A)}$ acts transitively on $\Delta_1$. To do this, 
observe that if $x\in S_1$ and $b\in A$ is positive and invertible,
then there exists an  invertible operator $R\in \Cal L_A(H_A)$ such that $R(x)=xb$: 
take $R=\theta_{xb , x} + 1-\theta_{x , x}$. It is straightforward to verify that 
$R$ is invertible and $R(x)=xb$. Now take $x,y \in \Delta_1$. Then $x'=x<x,x>^{-1/2}$ 
and $y'=y<y,y>^{-1/2}$ lie in $S_1$. Since the unitary group of $\Cal L_A(H_A)$ 
acts transitively on $S_1$, there exists a unitary $U$ such that $U(x')=y'$. 
Let $R_1,R_2$ be invertibles in $\Cal L_A(H_A)$ such that $R_1(x')=x'<x,x>^{1/2}=x$
and $R_2(y')=y$. Then $R_2UR_1^{-1}(x)=y$. Finally, we can see that $G_{\Cal L_A(H_A)}$ 
acts transitively on any $\Delta_q$: if $x,y \in \Delta_q$, pick $x_0, y_0 \in S_1$ 
such that $x_0q=x$ and $y_0q=y$ (recall that $m_q$ is surjective). If $R$ is an 
invertible operator such that $R(x_0)=y_0$, then also $R(x)=y$.

\noi
The transitivity of the action of $G_{\Cal L_A(H_A)}$ on both $\Delta_1$ and $\Delta_p$ makes it sufficient to find local cross sections for $m_p$ around
a fixed element of $\Delta_p$, say $e_1p$. Using the action the cross section can be carried over any point in $\Delta_p$. Define the C$^\infty$ (affine) map
$$
\sigma: \Delta_p(H_A) \to H_A \ , \quad \sigma(x)=x+e_2(1-p),
$$
for $x=(a_n) \in \Delta_p(H_A)$.
Note that $m_p(\sigma (x))=x$ and
$$
<\sigma(x),\sigma(x)>= <x,x>+1-p+(1-p)a_2+a_2^*(1-p).
$$
Note also, that
$$
\align
\|(1-p)a_2\|^2 & \le \|a_2\|^2=\|pa_2^*a_2p\| \le \|p(1-a_1)^*(1-a_1)p +\sum_{n\ge 2}pa_n^*a_np\| \\
            &=\|x-e_1p\|.
\endalign
$$
Fix any scalar$0<\delta <1$ and let $\Cal V_\delta$ be the open
neighbourhood of $e_1p$
$$
\Cal V=\{x\in \Delta_p(H_A): <x,x> > \delta p \hbox{ and } \|x-e_1p\|<\delta/2\}.
$$
Then it is clear (using the inequality above) that if $x\in \Cal V_\delta$ then $\sigma (x) \in \Delta_1(H_A)$. In other words, $\sigma$ is a local cross section for $m_p$ as claimed. \rect

\medskip
\noi \smc \bf Remark 3.5 \rm \rm
Suppose that $A$ is unital.
The fibre $F$ of $m_p: \Delta_1(H_A)\to \Delta_p(H_A)$
(over $e_1p$) is the space
$$
F=\{x=(a_n)\in \Delta_1(H_A): a_1p=p \hbox{ and } a_np=0 \hbox{ for all } n\ge 2\}.$$
Note that therefore $F$  splits as 
$$
F=\{a=a_1\in A: ap=p\}\times \{x'=(a_n)_{n\ge 2}: a_np=0 \hbox{ for all } n\ge 2\}.
$$
The first factor is clearly a convex set, and the second one is diffeomorphic to
$\Delta_{1-p}(H_A)$. Indeed, if $x=(a_n) \in F$,
$$
<x,x>=\sum a_n^*a_n=p +\sum_{n\ge 2}(1-p)a_n^*a_n(1-p).
$$
Since $<x,x>$ is invertible in $A$, it follows that $\sum_{n\ge 2}(1-p)a_n^*a_n(1-p)$ 
must be invertible in $(1-p)A(1-p)$.
In particular, it follows that $F$ is homotopically equivalent to 
$\Delta_{1-p}(H_A)$.
\medskip

One can relate the homotopy groups of $\Cal P(H_A)$ to the
homotopy groups of $\Cal U_{pAp}$, and eventually compute them.
\medskip
\noi \smc \bf Corollary 3.6 \rm \sl
Fix $x_0 \in S_p (H_A)$ and let $[x_0]$ be the submodule of $H_A$ generated by
$x_0$. Then for all $n \ge 1$,
$$
\pi_n (\Cal P(H_A), [x_0]) \simeq \pi_{n-1} (\Cal U_{pAp},p).
$$
In particular, if $A$ is a von Neumann algebra and $p $ is a properly infinite 
projection of $A$, then the connected component of $[x_0]$ in $\Cal P(H_A)$ is 
contractible. 
\smallskip \noi Proof. \rm
The result follows easily from the homotopy exact
sequence of the fibre bundle $\rho :S_p(H_A) \to \Cal P(H_A)_{[x_0]}$,
with fibre $\Cal U_{pAp}$.

If  $A$ is a von Neumann algebra and $p$ is a properly infinite projection, then [7] 
$\Cal U_{pAp}$ is contractible. Therefore the connected component of $[x_0]$  
in $\Cal P(H_A)$ has trivial homotopy groups. Since it is homeomorphic to 
$\Cal E_{ x_0}$, which is a differentiable manifold modeled on a Banach 
space (namely, $\Cal L_A(H_A)$), it follows [17] that it is contractible. \rect
\medskip
\noi \smc \bf Remark 3.7 \rm  \rm
In [25] Zhang  proved that if $A$ is a non elementary
\csta \ ($\ne \Cal K$ or $M_n(\zC)$) with real rank zero and stable
rank one, and $p$ is any projection in $A$,
then for $k\ge 1$
$$
\pi_{2k}(\Cal U_{pAp})\simeq K_1 (A)\simeq \pi_{2k+1}(\Cal P(A))$$
and
$$
\pi_{2k+1}(\Cal U_{pAp})\simeq K_0 (A)\simeq \pi_{2k+2}(\Cal P(A)).$$
Here $\Cal P(A)$ denotes the space of selfadjoint projections of $A$,
which agrees with the projective space of the module $X=A$.
These results, applied for both $A$ and $A\otimes \Cal K$ imply 3.6
for this class of algebras.
\medskip

If $A=B(H)$ with $H$ an infinite dimensional Hilbert space, then the connected 
components of  $\Cal P(H_{B(H)})$ corresponding to infinite rank projections are 
contractible. Note that because infinite rank projections in $B(H)$ are
equivalent, a single sphere, $S_1 (H_A)$, suffices to cover the whole connected component. 
The other connected components of $\Cal P(H_{B(H)})$, corresponding to 
finite rank projections $p$, can be parametrized by  the ranks of the 
projections.  If rk$(p)=n$, $1\le n <\infty$, then $\Cal P(H_{B(H)})_n$ has
trivial $\pi_0$ and $\pi_1$ and $\pi_2$ equal to $\zZ_n$.

More generally, if $A$ is a von Neumann algebra, $\pi_0(\Cal P(H_A))$ is 
parametrized by the equivalence classes of projections in $A$.

If $p$ is finite in $A$, then again $\pi_1(\Cal P(H_A),[x_0])$ ($x_0 \in S_p(H_A)$) 
is trivial and $\pi_2$ depends on the type decomposition of the finite algebra $pAp$. 
Suppose that $pAp= q_c(pAp) \oplus_{n\in \zN} q_n (pAp)$, with $q_c$ and $q_n$ the 
central projections of $pAp$ decomposing it in its type II$_1$ and I$_n$ parts, 
one has
$$
\pi_2(\Cal P(H_A),e_1p)\simeq C(\Omega_c, \zR) \oplus_{n\in \zN} C(\Omega_n,\zZ),
$$
where $\Omega_c$ and $\Omega_n$ denote the Stone spaces of the centres of
$q_c(pAp)$ and $q_n(pAp)$, respectively. 
\medskip

\bigskip
\noi \bf 4. Partial isometries of $\Cal L_A(H_A)$.  \rm
\medskip
If $x_0\in S_p(X)$, then $S_p(X)$ is (isometrically) isomorphic 
to $S_{\theta_{x_0, x_0}}(\Cal L_A(X))$. Indeed, any partial isometry $V$ in 
$\Cal L_A(X)$ with initial space $\theta_{x_0, x_0}$ is
of the form $V=\theta_{V(x_0), x_0}$. Then
$$
S_{\theta_{x_0, x_0}}(\Cal L_A(X)) \ni V \mapsto V(x_0) \in S_p(X)
$$
has inverse $y \mapsto \theta_{y , x_0}$. These maps are clearly isometric
since $\| \theta_{y, x_0} -\theta_{y', x_0}\|=\|\theta_{y-y', x_0}\|=\|y-y'\|$, where the last equality holds
because $y-y'=(y-y')p=(y-y')<x_0,x_0>$.

Let us consider from now on the case $X=H_A$. The space $S_1(\Cal L_A(H_A))$ 
has particular interest because if consists of all 
isometries of $H_A$. Recall the  bundle
$$
\rho :S_1(\Cal L_A(H_A))\to \Cal E_1$$
from the space of isometries onto the space of projections of $\Cal L_A(H_A)$ which are
equivalent to the identity, namely $\rho (V)= VV^*$. The fibre of this map is the 
unitary group of $\Cal L_A(H_A)$. It follows that $\pi_0 (S_1(\Cal L_A(H_A)))=
\pi_0 (\Cal E_1)$.
Now $\Cal E_1$ can be thought of as the set of (orthogonally) complemented 
submodules of $H_A$ which are isomorphic to $H_A$. Recall Kasparov's 
stability theorem, which states that any countable generated Hilbert module is 
isomorphic to a submodule of $H_A$ with complement isomorphic to the full $H_A$ [13]. 
Since the unitary group of $\Cal L_A(H_A)$ is connected, unitary equivalence 
classes of such submodules correspond to connected components of $\Cal E_1$. 
Therefore one has a natural bijection
$$
\{ \hbox{ classes of countably generated Hilbert modules over } A\} \leftrightarrow \pi_0(S_1(\Cal L_A(H_A))).
$$
Note that $\pi_0(S_1(\Cal L_A(H_A))$ has a natural semigroup structure, which can be induced on the set
of classes of countably generated $A$-modules.
\medskip
\noi \smc \bf Proposition 4.1 \rm \sl
Suppose that $P\in \Cal K_A(H_A)$ is a (compact) projection. Then the map
$$
m_P :S_1(\Cal L_A(H_A))\to S_P(\Cal L_A(H_A))
$$
is surjective, and therefore a C$^\infty$  fibre bundle.
\smallskip \noi Proof. \rm
Recall that $\Cal K_A(H_A)=\Cal K \otimes A =\lim_n M_n(A)$.
Since $P$ is compact, it is unitarily equivalent to a matrix projection $P_0$
in $M_n(A)$ for some $n$. It clearly suffices to show that $m_{P_0}$ is a
fibre bundle. Therefore we may suppose $P$ a matrix. Let us see first that
$m_P$ is onto. Pick $V\in S_P(\Cal L_A(H_A))$. Let $Q=VV^*$ be the final
projection of $V$. Again, since $Q$ is compact it is unitarily equivalent
to a matrix projection, $Q_0=UQU^* \in M_n(A)$ (without loss of generality
we may suppose the same size for $P$ and $Q_0$). Then the partial isometry
$UV$ with initial projection $P$ and final projection $Q_0$ is also an
$n\times n$ matrix. Let $W$ be the operator with
$$
\left( \matrix UV & 1-Q_0 \\ 1-P & V^*U^*  \endmatrix \right)
$$
on the first $2n\times 2n$ corner, and the identity on the rest. Then $W$
is a unitary operator in $\Cal L_A(H_A)$, satisfying $WP=UVP=UV$. Therefore
$V=U^*WP=m_P(U^*W)$, with $U^*W \in \Cal U_{\Cal L_A(H_A)} \subset S_1
(\Cal L_A(H_A))$. \rect
\medskip
\noi \smc \bf Remark 4.2 \rm  \rm
In 4.1 it is shown in fact that $m_P(\Cal U_{\Cal L_A(H_A)})=
S_P(\Cal L_A(H_A))$, and therefore that the restriction of $m_P$ to the unitary
group of $\Cal L_A(H_A)$ (= the connected component of the identity in
$S_1(\Cal L_A(H_A))$) is a fibre bundle. In particular, this implies that $S_P(\Cal
L_A(H_A))$ is connected. Moreover
\medskip
\noi \smc \bf Corollary 4.3 \rm \sl
Suppose that $P\in \Cal K\otimes A$ is a projection, then the space
$S_P(\Cal K\otimes A)$ of partial isometries in $\Cal K \otimes A$ with
initial space $P$ is contractible.
\smallskip \noi Proof. \rm
Note that $S_P(\Cal K\otimes A)$ coincides with $S_P(\Cal L_A(H_A))$. Consider
the restriction of $m_p$ to the unitary group of $\Cal L_A(H_A)=M(\Cal
K\otimes A)$,
$$
m_P :\Cal U_{M(\Cal K\otimes A)} \to S_P(\Cal K\otimes A).$$
This map is a bundle with fibre equal to $\{V:VP=P\}$. Clearly this set
identifies with the unitary group of the submodule $R(P)^\perp \subset H_A$.
Since  $P$ is compact, by [16,1.10] there exists an isometry $W$ in $M(A\otimes \Cal K)$ such that $1-P=WW^*$. Then $R(P)^\perp
\simeq H_A$, and therefore the fibre is contractible.
Then $S_P(\Cal K\otimes A)$ is a differentiable manifold, and the
base space of a contractible fibre bundle with contractible fibre. It
follows from the already cited result of [17] that $S_P(\Cal K \otimes A)$ is contractible. \rect
\medskip
As a consequence, putting $P=\theta_{x, x}$ for some $x\in H_A$ with $<x,x>=p$, one obtains the proof of the contractability of $S_p(H_A)$, because, as noted at the beginning of this section, $S_{\theta_{x, x}}(A\otimes \Cal K) \simeq S_p(H_A)$.
\medskip
In [16] Mingo defines an index map (generalizing the index for Fredholm
operators in Hilbert spaces) which induces isomorphism
$$
\{ \hbox{ Fredholm partial isometries of } H_A \} \to K_0 (A)$$
where Fredholm operators of $H_A$ are elements of $\Cal L_A(H_A)$ which have closed complemented range, with finitely generated kernel and cokernel.
If one restricts this  map to the semigroup of classes (Fredholm)isometries, one obtains  $-K_0^{+}(A)$. In other words, (classes of) Fredholm isometries correspond to finitely generated $A$-modules. The remaining part of $S_1$, which could be called the semi-Fredholm isometries, correspond with the infinite (countable) generated $A$-modules.

\bigskip
\noi \bf 5. Projective space of a selfdual module.  \rm
\medskip
In this section we consider the case when $A$ is a von Neumann algebra and
$X$ is selfdual [18]. Then $\Cal L_A(X)$ is a von Neumann algebra with the 
same centre as $A$.

Let us first state the following result, which can be proved analogously as in 6.2 of [1].
\medskip
\noi \smc \bf Lemma 5.1 \rm \sl
Let $M$ be a von Neumann algebra of type II$_1$ (resp. I$_n$, $n<\infty$), $Tr$ its center valued trace and $e \in M$ a projection. Let $i:\Cal U_{eMe} \hookrightarrow \Cal U_M$ be the inclusion map $i(eue)=eue+1-e$. Then the image $i^*(\pi_1(\Cal U_{eMe}))\subset \pi_1 (\Cal U_M)$ is the group of multiples of $Tr(e)\hat{ }$, where, as in [11], $\pi_1(\Cal U_M)$
is identified with the additive group $C(\Omega,\zR)$ (resp. $C(\Omega, \zZ)$),
$\Omega$ is the Stone space of the centre of $M$ and $Tr(e)\hat{ }$
is the Gelfand transform of $Tr(e)$. \rect
\rm
\medskip
Denote by $\Cal Z(A)$ the centre of $A$, and $\Cal Z(A)_{sa}$ the subspace of selfadjoint elements of $\Cal Z(A)$
\medskip
\noi \smc \bf Remark 5.2 \rm \sl \rm
In [1] we computed the fundamental group of the sphere $S_1(X)$.
The procedure to compute the fundamental group of $S_p(X)$ for any projection $p$ is similar. Fix an element $x_0 \in S_p(X)$ and consider the bundle
$$
\pi_{x_0} :\Cal U_{\Cal L_A(X)} \to S_p(X) , \ \pi_{x_0}(U)=U(x_0)
$$
with fibre
$$
F=\{V\in \Cal U_{\Cal L_A(X)}: Vx_0=x_0\}
$$
introduced in 2.8. $\Cal L_A(X)$ and $A$ have the same central projections decomposing
them in their type I,II and III parts (see [19]), although it clearly can happen $A$
finite with $\Cal L_A(X)$ infinite. Using these projections the spheres $S_p(X)$ split and
one is set in the case when $\Cal L_A(X)$ and $A$ are of one and the same definite type. 
\medskip
If $A$
is properly infinite, one has the following:
\medskip
\noi \smc \bf Theorem 5.3 \rm \sl
If $A$ is properly infinite, then the connected components of $S_p(X)$ are contractible.
\smallskip \noi Proof. \rm
If $A$ is properly infinite, then $\Cal L_A(X)$ is properly infinite. In [1] 
it was proven that $A$ infinite implies $\Cal L_A(X)$ infinite. Suppose that 
$\Cal L_A(X)$ is not properly infinite, then there would exist a finite central 
projection in $\Cal L_A(X)$, which would imply the existence of a finite central 
projection in $A$. Therefore in the bundle $\pi_p$ one has that $\Cal U_{\Cal L_A(X)}$ 
is contractible ([7]). The fibre $F$ identifies with the unitary group of the 
submodule $Y=[x_0]^{\perp}$, which is also a selfdual module over $A$. If $Y$ 
is non trivial, then again $F$ is contractible, and the connected component of 
$[x_0]$ in $S_p(X)$ is contractible, being a manifold with trivial homotopy groups. If $Y$ is trivial, then $X=[x_0]$, which implies that $X$ is isomorphic to $pAp$ (as $pAp$ modules). Then $F$ reduces to the trivial group, and the result follows  (in this case, $S_p(X)$ is just the isometries of $pAp$, whose connected components are homeomorphic to the unitary group of $pAp$). \rect
\medskip
If $A$ is finite, then $S_p(X)$ is connected. Indeed, the projections $\theta_{x, x}$ for
$x\in S_p(X)$ are equivalent and, in this case, finite. Then they are unitarily equivalent,
which implies that the action of $\Cal U_{\Cal L_A(X)}$ is transitive in $S_p(X)$ (see [1]).
In this case $\Cal L_A(X)$ can be either finite or infinite, and 
there exists a central projection $p_f$ in $A$ such that $p_f \Cal
L_X(X)$ is finite and $(1-p_f) \Cal L_A(X)$ is properly infinite. 
Note that the first algebra identifies with  $\Cal
L_{Ap_f}(Xp_f)$, and the second with $\Cal L_{A(1-p_f)}(Xp_f)$.
The spheres and the projective space split, $S_p(X)\cong S_{pp_f}(Xpf)\times S_{p(1-p_f)}(X(1-p_f))$ and 
$\Cal P(X)\cong \Cal P(Xp_f)\times
\Cal P(X(1-p_f))$. So one may consider separately the cases when
$\Cal L_A(X)$ is finite or properly infinite.

If $\Cal L_A(X)$ is properly infinite, 
the unitary group of $\Cal L_A(X)$ is contractible. Moreover, since $\theta_{x, x}$ is finite, 
$(1-\theta_{x, x})\Cal L_A(X)(1-\theta_{x, x})$ is also properly infinite. The unitary group
of this algebra identifies with the fibre of the bundle $\pi_x$ over the sphere $S_p(X)$, 
$\{V\in \Cal U_{\Cal L_A(X)}: V(x)=x\}$.
Therefore one has the following:
\medskip
\noi \smc \bf Proposition 5.4 \rm \sl
If $A$ is finite, and $p_f$ is defined as above, then the sphere $S_p(X)$ is homotopically 
equivalent to its finite part $S_{pp_f}(Xp_f)$.
\smallskip \noi Proof. \rm
The proof proceeds as in the above result, observing that the infinite parts corresponding to the central
projection $1-p_f$ are contractible. \rect

The fundamental groups of the spheres in the finite case were computed in [1].
If $\Cal L_A(X)$ is of type II$_1$, $S_p(X)$ is connected and one has the tail of the homotopy exact sequence
$$
\pi_1(F)\hookrightarrow^{i^*} \pi_1(\Cal U_{\Cal L_A(X)})\to \pi_1(S_p(X)) \to 0,
$$
where $F=\{V\in \Cal U_{\Cal L_A(X)}: Vx_0=x_0\}$.
Applying 5.1 above and the results on [11] characterizing the fundamental groups of the unitary groups of von Neumann algebras, one obtains that the image of $i^*$ equals the (additive)
group $\{z(1-Tr(\theta_{x_0, x_0})): z\in \Cal Z(A)_{sa}\}$, and therefore
$$
\pi_1 (S_p(X))=\Cal Z(A)_{sa}/\{z(1-Tr(\theta_{x_0, x_0})): z\in \Cal Z(A)_{sa}\},
$$
for any $x_0 \in S_p(X)$.

Proceeding analogously, if $\Cal L_A(X)$ is of type I$_n$ ($n<\infty$), one
has
$$
\pi_1(S_p(X))=C(\Omega,\zZ)/\{f(1-Tr(\theta_{x_0,x_0})\hat{ }): f\in C(\Omega,\zZ)\},
$$
where $\Omega$ is the Stone space of $\Cal Z(A)$. Note that also in this case $S_p(X)$ is connected.
\medskip

We will show that the connected components of the
projective space $\Cal P(X)$ have trivial fundamental group. This will be done again using the results of Handelmann [11], and the principal bundle
$$
\Cal U_{\Cal L_A(X)} \to \Cal P(X)_{[x]} \ ,  \ U\mapsto [U(x)]
$$
for a given $[x] \in \Cal P(X)$.
In fact, this result will follow from considering the general case, of a
von Neumann algebra $M$ and an arbitrary projection $p\in M$.
\medskip
\noi \smc \bf Theorem 5.5 \rm \sl
The unitary orbit $\Cal U_M(p)=\{upu^*: u\in \Cal U_M\}$ is simply connected.
\smallskip \noi Proof. \rm
Clearly $\Cal U_M(p)$ is connected.
Consider the principal bundle $\pi_p$
$$
\pi_p :\Cal U_M \to \Cal U_M (p) \ , \ \pi_p(u)=upu^*,
$$
with fibre $\{ v\in \Cal U_M : vpv^*=p\}=\Cal U_M \cap \{p\}'$. Note that
the fibre is the unitary group of the commutant $M\cap \{p\}'$, and is therefore connected.

If $q$ is a projection in the center of $M$, then the unitary orbit factors
as the unitary orbit of $qp$ under the action of the unitary group of $qM$,
times the unitary orbit of $(1-q)p$ under the action of the unitary group of
$(1-q)M$. Therefore, using the type decomposition central projections of $M$,
one may consider separately the cases $M$ properly infinite, type II$_1$ and
type I$_n$, for $n<\infty$.

If $M$ is properly infinite, $\Cal U_M$
has trivial $\pi_1$-group. It follows that in this case also the unitary
orbit has trivial $\pi_1$-group.

Suppose now that $M$ is either of type II$_1$ or type I$_n$ with
$n<\infty$. Note that the fibre $\Cal U_M\cap\{p\}'$ factors as
$\Cal U_{pMp} \times \Cal U_{(1-p)M(1-p)}$. Let us show that the
homomorphism $i^*:\pi_1(\Cal U_M\cap \{p\}') \to \pi_1(\Cal U_M)$
induced by the inclusion map $i :\Cal U_M\cap \{p\}' \hookrightarrow \Cal U_M$
is surjective. The image of $i^*$ contains both
$i^*(\pi_1 (\Cal U_{pMp}))$ and $i^*(\pi_1 (\Cal U_{(1-p)M(1-p)}))$, which,
by the above lemma, are generated by the multiples of $Tr(p)\hat{ }$ and $1-Tr(p)\hat{ }$ respectively. Therefore $i^*$ is surjective. Using the homotopy exact sequence of the bundle $\pi_p$,
$$
\dots \pi_1 (\Cal U_M\cap\{p\}')\to^{i^*} \pi_1 (\Cal U_M) \to \pi_1 (\Cal U_M(p)) \to \pi_0 (\Cal U_M\cap\{p\}')=0,
$$
since $i^*$ is surjective, it follows that $\pi_1(\Cal U_M(p))$ is trivial. \rect
\medskip
\noi \smc \bf Corollary 5.6 \rm \sl
If $M$ is a properly infinite von Neumann algebra and $p\in M$ is a properly infinite projection with $1-p$ also properly infinite, then $\Cal U_M(p)$ is contractible.
\smallskip \noi Proof. \rm
In this case, one has that in the proof of the preceeding result, the structure 
group  $\Cal U_{pMp} \times \Cal U_{(1-p)M(1-p)}$ and the space of the bundle, 
$\Cal U_M$, are both contractible. The proof follows because $\Cal U_M(p)$ 
is a differentiable manifold. \rect
\medskip
\noi \smc \bf Corollary 5.7 \rm \sl
Let $X$ be a selfdual right C$^*$-module over the von Neumann algebra $A$.
Pick $[x]\in \Cal P(X)$, with $x\in S_p(X)$, then the connected component
$\Cal P(X)_{[x]}$ of $[x]$ has trivial $\pi_1$-group.

If  $A$ is properly infinite, then
$$
\pi_n(\Cal P(X),[x])\simeq \pi_{n-1}(\Cal U_{pAp},p).
$$
If moreover $p$ is properly infinite, then $\Cal P(X)_{[x]}$ is contractible.
\smallskip \noi Proof. \rm The proof follows by observing that $\Cal P(X)_{[x]}$ is
homeomorphic to the unitary orbit $\Cal U_{\Cal L_A(X)}(\theta_{x, x})$. If $A$ is properly
infinite, then $S_p(X)$ is contractible, and the second statement follows considering the
bundle $\rho$. Finally, if $p$ is properly infinite, $\Cal U_{pAp}$ is contractible,
and therefore $\Cal P(X)_{[x]}$ is contractible.  \rect

If $A$ is of type II$_1$, then the second homotopy group of the projective
space is non trivial. We shall consider this fact next.
Recall that one can consider separately the cases when $\Cal L_A(X)$ is
finite and properly infinite, using the central projection $p_f$
defined before.
First, we need the following result:
\medskip
\noi \smc \bf Lemma 5.8 \rm \sl
Suppose that $B$ is a von Neumann algebra of type II$_1$ and $p\in B$ a projection. Then the inclusion
$i:\Cal U_{pBp} \to \Cal U_B$, $i(pup)=pup+1-p$ induces the (additive) group homomorphism between the $\pi_1$ groups
$$
i^* :\pi_1(\Cal U_{pBp},p) \to \pi_1(\Cal U_B,1) , \ i^*(xp)=Tr(xp)
$$
where $Tr$ is the centre valued trace of $B$, and $\pi_1(\Cal U_B)$ (resp. $\pi_1(\Cal U_{pBp})$) is identified with
$\Cal Z(B)_{sa}$ (resp. $\Cal Z(pBp)_{sa}=\Cal Z(B)_{sa}p$).
\smallskip \noi Proof. \rm
In [11] and [4] it was shown that the classes of the loops $e^{itq}$, with $q$ a projection, generate the fundamental group of the unitary group. Moreover, the isomorphism identifying $\pi_1(\Cal U_B)$
with $\Cal Z(B)_{sa}$ takes the class of the loop $e^{itq}$ to the element $Tr(q)$. The analogous fact holds for $pBp$, except that
one considers projections $q \le p$, and one uses $Tr_{pBp}$ the centre valued trace of $pBp$. Now, it is clear that $Tr_{pBp}(pxp)=Tr(pxp)p$. Therefore, if $q\le p$,
$Tr_{pBp}(q)=Tr(q)p \in \Cal Z(B)p$ which corresponds to the class of the loop $e^{itq}p$ in $\Cal U_{pBp}$, is  mapped to the class of the loop $e^{itq}p +1-p=e^{itq}$ in $\Cal U_B$. This class corresponds
to the element $Tr(q)$ in $\Cal Z(B)_{sa}$. That is, $i^*(xp)=Tr(xp)$ for all $x\in \Cal Z(B)_{sa}p$.\rect

\medskip
\noi \smc \bf Proposition 5.9 \rm \sl
Suppose that $A$ is of type II$_1$. Fix $[x] \in \Cal P(X)$ and denote
by $e$ the projection $\theta_{x, x}$.
\item{ a)} If $\Cal L_A(X)$ is properly infinite, then
$$
\pi_2(\Cal P(X),[x])\simeq \Cal Z(A)_{sa}e.
$$
\item{ b)} If $\Cal L_A(X)$ is finite, denote by $Tr$ the centre
valued trace of $\Cal L_A(X)$. Then
$$
\pi_2(\Cal P(X),[x])\simeq \{(a,b)\in \Cal Z(A)_{sa}e\times \Cal Z(A)_{sa}(1-e): Tr(a+b)=0\}.
$$
Both sets on the right hand are considered as additive groups.
\smallskip \noi Proof. \rm
Recall the homotopy exact sequence of the bundle $\Cal U_{\Cal
L_A(X)} \to \Cal P(X)_{[x]}$, with fibre equal to the unitary group
of $\Cal L_A(X) \cap \{e\}'$. The fibre is homeomorphic to the
product of the unitary groups of  $e\Cal L_A(X) e$ and $(1-e)\Cal
L_A(X) (1-e)$. In case a),
since $e$ is finite and $\Cal L_A(X)$ is properly infinite,  
the unitary group of $(1-e)\Cal L_A(X) (1-e)$ is contractible. Therefore in
$$
\pi_2(\Cal
L_A(X),1) \to \pi_2(\Cal P(X)_{[x]},[x])  \to \pi_1(\Cal
U_A,1)\times \pi_1(\Cal U_{(1-e)\Cal L_A(X)(1-e)}, 1-e) \to
\pi_1(\Cal L_A(X),1)
$$
one has that $\pi_i (\Cal U_{\Cal L_A(X)})$, $i\ge 0$
and $\pi_1(\Cal U_{(1-e)\Cal L_A(X)(1-e)})$ are trivial. 

Then
$\pi_2(\Cal P(X)_{[x]},[x])\simeq \pi_1(\Cal U_A,1)\simeq \Cal
Z(A)_{sa}$.

In case b), i.e. $\Cal L_A(X)$ finite, Schr\"{o}der [23] proved
that $\pi_2(\Cal L_A(X))=0$. In this case $\pi_1(\Cal U_{(1-e)\Cal
L_A(X)(1-e)})\simeq \Cal Z(A)_{sa}(1-e)$, i.e. the selfadjoint
elements of the centre of $(1-e)\Cal L_A(X)(1-e)$. Therefore one
has $$ \matrix
      &                            & \partial &                                                & i^* &                &       \\
0 \to & \pi_2(\Cal P(X)_{[x]},[x]) & \to & \Cal Z(A)_{sa}e\times \Cal Z(A)_{sa}(1-e) & \to & \Cal Z(A)_{sa} & \to 0.
\endmatrix
$$ 
By 5.8, the morphism $i^*$ is
given by $i^*(a,b)=Tr(ae+b(1-e))$. On the other hand, the sequence
above shows that the map $\partial :\pi_2(\Cal P(X)_{[x]},[x]) \to
\Cal Z(A)_{sa}\times \Cal Z(A)_{sa}(1-e)$ is injective. Therefore
$$ \pi_2(\Cal P(X)_{[x]},[x])\simeq \partial (\pi_2(\Cal
P(X)_{[x]},[x]))=\ker i^* , $$ which ends the proof. \rect
\medskip
 If $A$ is a factor, then $\Cal L_A(X)$ is either finite
or properly infinite. In [1] it was noted that $\Cal L_A(X)$ is
finite if and only if $X$ is finitely generated. In both
situations, it follows from the preceeding result that the second
homotopy group of $\Cal P(X)$ is isomorphic to $(\zR, +)$.

\medskip

\noi \smc \bf Remark 5.10  \rm \sl \rm Let us recall the example
where $X$ is a von Neumann algebra $B$ containing $A$, and there
exists a finite index conditional expectation $E:B\to A$, which
will be automatically normal. Then ([5]) $B$ is a selfdual
$A$-module and the results above apply. The algebra $\Cal L_A(B)$
is isomorphic to the Jones extension $B_1$, namely the von Neumann
algebra generated by $B$ and the Jones projection $e$ (which is
the map $E$ regarded as an element on $\Cal L_A(B)$). As shown in
[1], if $A\subset B$ are of type II$_1$, $\pi_1(S_1^E(B))$ is
isomorphic to the additive group $$ Z(A)_{sa}/\{xTr(1-e): x\in
x\in Z(A)_{sa}\} $$ where $Z(A)_{sa}$ denotes the set  of
selfadjoint elements of the centre $Z(A)$ of $A$. Note that
$Tr(e)$ equals the inverse of $E([E])$, where $[E]$ is the center
valued index (in the centre of $B$) of the expectation $E$ (see
[5]).

\medskip

\noi \smc \bf Remark 5.11  \rm \sl \rm Suppose now that $A\subset
B$ is a finite index inclusion of factors of type II$_1$. As in
the remark above, regard $B$ as a (selfdual) C$^*$-module over
$A$. In this case the first homotopy group of $\Cal P_A(B)$ is
trivial, and the second is $\zR$. One recovers the index of the
inclusion as the generator of $\partial (\pi_2(\Cal P_A(B)))
\subset \zR^2$. Indeed, in the proof of 5.9 it was shown that
$\partial (\pi_2(\Cal P_A(B)))=\{(s,t)\in \zR^2:
Tr(se+t(1-e))=0\}$. Here $Tr(e)=[B:A]^{-1}$, where $[B:A]$ is the
Jones index of the inclusion. A generator for $\partial
(\pi_2(\Cal P_A(B)))$ is then the pair $(1, 1-[B:A])$
\medskip

\noi \bf References. \rm
\bigskip
\item { [1] } E. Andruchow, G. Corach and D. Stojanoff; Geometry of
the sphere of a Hilbert module, Math. Proc. Cambridge Phil. Soc. (to appear)
\item { [2] } E. Andruchow, G. Corach and D. Stojanoff; Projective space of a 
C$^*$-algebra, Integral Equations and Operator Theory (to appear).
\item { [3] } E. Andruchow, A. R. Larotonda, L. Recht and D. Stojanoff;
Infinite dimensional homogeneous reductive spaces and finite index
conditional expectations, Illinois Math. J. 41 (1997), 54-76.
\item { [4] } H.Araki, M.Smith and L. Smith; On the 
homotopical significance of the type of von Neumann
algebra factors, Commun. Math. Phys. 22 (1971),71-88.
\item { [5] } M. Baillet, Y. Denizeau and J.F. Havet; Indice d'une 
esperance conditionelle,  Comp. Math. 66 (1988), 199-236.
\item{ [6] } M. Breuer, A generalization of Kuiper's theorem to factors of type II$_\infty$, J. Math. Mech. 16 (1967), 917-925.
\item{ [7] } J. Br\"{u}ning, W. Willgerodt; Eine Verallgemeinerung
eines Satzes von N. Kuiper, Math. Ann. 220 (1976), 47-58.
\item { [8] } G. Corach, H. Porta and L. Recht; The geometry of spaces of
projections in \csta s, Adv. Math. Vol. 101 (1993), 59-77.
\item{ [9] } J. Cuntz and N. Higson, Kuiper's theorem for Hilbert modules;
Contemporary Mathematics 62 (1987), 429-435.
\item{ [10] } M. Frank and E. Kirchberg; Conditional expectations of finite index;
J. Oper. Theory 40 (1998) 87-111.
\item{ [11] } D. E. Handelman; K$_0$ of von Neumann algebras and
AFC$^*$-algebras, Quart. J. Math. Oxford (2) 29 (1978), 429-441.
\item { [12] } I. Kaplansky; Modules over operator algebras, Amer. J. Math. 
75 (1953), 839-858.
\item { [13] } G.G. Kasparov; Hilbert C$^*$-modules: theorems of Stinespring and
Voiculescu, J. Operator Theory 4 (1980), 133-150.
\item { [14] } E.C. Lance; Hilbert C$^*$-modules - a toolkit for operator 
algebraists; London Math. Soc. Lecture Notes Series 210, Cambridge 
University Press, Cambridge, 1995.
\item { [15] } A.R. Larotonda; Notas sobre variedades diferenciables, Notas
de Geometr\'\i a y topolog\'\i a 1, INMABB-CONICET, Universidad Nacional del
Sur, Bah\'\i a Blanca, Argentina, 1980.
\item { [16] } J.A. Mingo; $K$-theory and multipliers of stable $C^*$-algebras,
Trans. Amer. Math. Soc. vol. 299, (1987), 397-411.
\item{ [17]  } R.S. Palais; Homotopy theory of infinite dimensional manifolds, Topology 5 (1966), 1-16.
\item { [18] } W.L. Paschke; Inner product modules over B$^*$-algebras; 
Trans. Amer. Math. Soc. 182 (1973),443-468.
\item { [19] } W.L. Paschke; Inner product modules arising from compact 
automorphism groups of von Neumann algebras, Trans. Amer. Math. Soc. 
224 (1976), 87-102.
\item { [20] } S. Popa; Classification of subfactors and their endomorphisms, CBMS 86,
AMS (1995).
\item { [21] } H. Porta and L. Recht, Minimality of geodesics in Grassmann manifolds,
Proc. Amer. Math. Soc. 100 (1987), 464-466.
\item{ [22] } M.A. Rieffel; Induced representations of C$^*$-algebras, Adv. Math. 13 (1974), 176-257.
\item{ [23] } H. Schr\"{o}der; On the homotopy type of the regular group of a
W$^*$-algebra, Math. Ann. 267 (1984), 694-705.
\item{ [24] } D.R. Wilkins; The Grassmann manifold of a C$^*$-algebra, Proc. Royal Irish Acad. 90A (1990), 99-116.
\item { [25] } S. Zhang; Matricial structure and homotopy type of simple \csta s with
real rank zero, J. Operator Theory 26 (1991), 283-312.
\vglue.5truecm

\noindent{Esteban Andruchow}

\noindent{Instituto de Ciencias, UNGS,  San Miguel, Argentina}

\noindent  J.A. Roca 850, 1663 San Miguel \ \ Argentina 

\noindent{e-mail : eandruch\@ungs.edu.ar}

\bigskip

\noindent{Gustavo Corach}

\noindent{Depto. de Matem\'atica, FCEN-UBA, Buenos Aires, Argentina}

\noindent  Ciudad Universitaria, 1428 Buenos Aires  \ \  Argentina 

\noi and 

\noindent{Instituto Argentino de Matem\'atica, Buenos Aires,  Argentina}

\noindent  Saavedra 15, 1083 Buenos Aires \ \  Argentina

\noindent{e-mail: gcorach\@mate.dm.uba.ar}

\smallskip

\noindent{Demetrio Stojanoff}

\noindent{Depto. de Matem\'atica, FCE-UNLP,  La Plata,  Argentina}

\noindent 1 y 50, 1900 La Plata \ \  Argentina

\noindent{e-mail: demetrio\@mate.dm.uba.ar}

\bigskip
 
\noindent{Mathematical Subject Classification 1991: Primary 46L05, 58B20.}

\end

\centerline {\bf Esteban Andruchow\rm}
\centerline {Instituto de Ciencias, Universidad Nacional de General Sarmiento }
\centerline { J.A. Roca 850, 1663 San Miguel \ \ Argentina }
\medskip
\centerline {\bf  Gustavo Corach \rm}
\centerline {Instituto Argentino de Matem\'atica }
\centerline { Saavedra 15, 1083 Buenos Aires \ \  Argentina}
\centerline {and}
\centerline { Departamento de Matem\'atica, FCEN, Universidad de Buenos Aires }
\centerline { Ciudad Universitaria, 1428 Buenos Aires  \ \  Argentina }
\medskip
\centerline {\bf Demetrio Stojanoff\rm}
\centerline { Departamento de Matem\'atica, FCE, Universidad Nacional de La Plata }
\centerline{1 y 50, 1900 La Plata \ \  Argentina}

\end